\documentclass[reqno,10pt]{amsart}
\usepackage{amssymb}
\usepackage[usenames, dvipsnames]{color}
\usepackage{esint}
\usepackage{hyperref}
\usepackage{todonotes}
\usepackage{verbatim}
\usepackage{enumerate}

\usepackage{cite}
\usepackage[nobysame]{amsrefs}
\def\MR#1{}

\theoremstyle{plain}
\newtheorem{theorem}{Theorem}[section]

\theoremstyle{definition}

\theoremstyle{remark}

\newcommand{\supp}{\operatorname{supp}}
\newcommand{\dist}{\operatorname{dist}}

\numberwithin{equation}{section}

\newcommand{\bC}{\mathbb{C}}
\newcommand{\bN}{\mathbb{N}}

\newcommand{\bR}{\mathbb{R}}

\newcommand\cB{\mathcal{B}}

\newcommand\cD{\mathcal{D}}

\newcommand\cG{\mathcal{G}}

\makeatletter
\def\dashint{\operatorname%
{\,\,\text{\bf--}\kern-.98em\DOTSI\intop\ilimits@\!\!}}
\makeatother

\begin{document}
\title[Optimal estimates for transmission problems]{Optimal estimates for transmission problems including relative conductivities with different signs}

\author[H. Dong]{Hongjie Dong}
\author[Z. Yang]{Zhuolun Yang}
\address[H. Dong]{Division of Applied Mathematics, Brown University, 182 George Street, Providence, RI 02912, USA}
\email{Hongjie\_Dong@brown.edu}

\address[Z. Yang]{Institute for Computational and Experimental Research in Mathematics, Brown University, 121 South Main Street, Providence, RI 02903, USA}
\email{zhuolun\_yang@brown.edu}

\thanks{H. Dong is partially supported by Simons Fellows Award 007638, NSF Grant DMS-2055244, and the Charles Simonyi Endowment at the Institute for Advanced Study.}
\thanks{Z. Yang is partially supported by Simons Foundation Institute Grant Award ID 507536.}
%
%\keywords{parabolic equation, time fractional derivative, mean oscillation estimates, measurable coefficients}
%
%

\begin{abstract}
We study the gradient and higher order derivative estimates for the transmission problem in the presence of closely located inclusions. We show that in two dimensions, when relative conductivities of circular inclusions have different signs, the gradient and higher order derivatives are bounded independent of $\varepsilon$, the distance between the inclusions. We also show that for general smooth strictly convex inclusions, when one inclusion is an insulator and the other one is a perfect conductor, the derivatives of any order is bounded independent of $\varepsilon$ in any dimensions $n \ge 2$.
\end{abstract}

\maketitle

\section{Introduction and main result}

In this paper, we first study the second-order elliptic equations in divergence form with discontinuous coefficients in two dimensions
\begin{equation}\label{main}
L_{\varepsilon;r_1,r_2} u:= D_i(a(x) D_i u) = D_i f_i \quad \mbox{in}~~\cD,
\end{equation}
where $\cD$ is a bounded open subset of $\bR^2$,
$$
a(x) = k_1 \chi_{\cB_1} + k_2 \chi_{\cB_2} + k_0\chi_{\cB_0},
$$
$k_0=1$, $k_1,k_2,r_1,r_2 \in (0,\infty)$ are constants,
$$
\cB_1:= B_{r_1}(\varepsilon/2 + r_1,0), \quad \cB_2:= B_{r_2}(-\varepsilon/2 - r_2,0), \quad \cB_0:= \bR^2 \setminus (\overline{\cB_1 \cup \cB_2}),
$$
and $\chi$ is the indicator function.
The equation models the conductivity problem in composite material. The gradient of the voltage potential $u$ represents the electric field, and $a(x)$ is the conductivity which is a constant on each inclusion, and a different constant on the background matrix. It is significant from an engineering point of view to estimate the derivatives of the solutions.

In \cite{BASL}, Babu\v{s}ka et al. analyzed an analogous elliptic system, and numerically showed that, when the ellipticity constants are away from $0$ and infinity, the gradient of solutions remains bounded independent of $\varepsilon$, the distance between inclusions. When $\varepsilon = 0$, Bonnetier and Vogelius \cite{BV} proved that $|D u|$ is bounded for a fixed $k = k_1 = k_2$ away from $0$ and infinity. This result was extended by Li and Vogelius \cite{LV} to general second order elliptic equations in divergence form with piecewise H\"older coefficients and general shape of inclusions in any dimension. Furthermore, they established a stronger piecewise $C^{1,\alpha}$ control of $u$, which is independent of $\varepsilon$. Li and Nirenberg \cite{LN} further extended this global Lipschitz and piecewise $C^{1,\alpha}$ result to general second order elliptic systems in divergence form, including the linear system of elasticity. Some higher order derivative estimates in dimension $n=2$ were obtained in \cites{DZ,DL,JiKang}.

On the other hand, if $k_1,k_2$ are allowed to be $0$ or $\infty$, it was shown in \cites{Kel, BudCar, Mar} that the gradient of solutions generally becomes unbounded as $\varepsilon \to 0$. For the perfect conductivity problem ($k_1 = k_2 = \infty$), it has been proved that the generic blow-up rate of $|Du|$ is $\varepsilon^{-1/2}$ in two dimensions, $|\varepsilon \log \varepsilon|^{-1}$ in three dimensions, and $\varepsilon^{-1}$ in dimensions greater than three; see \cites{AKLLL,AKL,Y1,Y2,BLY1,BLY2}. These bounds were shown to be optimal and  are independent of the shape of inclusions, as long as the inclusions are relatively strictly convex. Moreover, more detailed characterizations of the singular behavior of $\nabla u$ have been obtained. For further works on the perfect conductivity problem and closely related works, see e.g. \cites{ACKLY,BT1,BT2,DL,AKKY,KLY1,KLY2,L,LLY,LWX,BLL,BLL2,DZ,KL,CY,ADY,Gor,LimYun,Kang} and the references therein.

For the insulated conductivity problem ($k_1 = k_2 = 0$), it was shown in \cites{AKLLL,AKL} that the optimal blow-up rate is $\varepsilon^{-1/2}$ in two dimensions. The proof uses a harmonic conjugate argument to link the insulated conductivity case to the perfect conductivity case, which fails in dimensions greater than two. For the higher dimensional case, Bao, Li, and Yin in \cite{BLY2} established an upper bound of order $\varepsilon^{-1/2}$. Yun in \cite{Y3} proved the optimal blow-up rate on the shortest line segment connecting two spherical inclusions in three dimensions is $\varepsilon^{\frac{\sqrt{2} - 2}{2}}$. Later, the upper bound $\varepsilon^{-1/2}$ was improved by Li and Yang in \cite{LY2} to be $\varepsilon^{-1/2 + \beta}$, for some $\beta > 0$. See also \cite{LY}. Weinkove in \cite{We} used a Bernstein-type argument to obtain a more explicit upper bound of $\beta$ in dimensions greater than three.
Dong, Li, and Yang in recent works \cites{DLY,DLY2} identified the optimal blow-up rate. They proved the optimal gradient estimate for a class of inclusions including balls and ``almost" optimal gradient estimate for general strictly convex inclusions. Unlike the perfect conductivity case, the optimal blow-up rate is related to the principal curvature of the inclusions.

Recently, among other results, Ji and Kang in \cite{JiKang} used some spectral properties of Neumann-Poincar\'e operators to study the problem \eqref{main} for the case when $0 < k_1 < 1$ and $k_2 > 1$, and proved that
$$
|D^m u| \le C \left( -\frac{(k_1+1)(k_2+1)}{(k_1-1)(k_2-1)}-1 + \sqrt{\frac{2(r_1 + r_2)\varepsilon}{r_1r_2}} \right)^{-m+1},\quad m=1,2,\ldots.
$$
In particular, when $k_1 \to 0$ and $k_2 \to \infty$, this implies
$$
|D^m u| \le C \varepsilon^{-(m-1)/2},\quad m=1,2,\ldots.
$$
In the first part of this paper, we apply the Green function method developed in \cite{DL} to show that $|D^m u|$ is in fact bounded independent of $\varepsilon$ for any $m=1,2,\ldots$. Compared to the method in \cite{JiKang}, we use $m$-th order finite differences to estimate the $m$-th order derivatives more precisely.

Let
$$
\alpha = \frac{k_1-1}{k_1 + 1}, \quad \beta = \frac{k_2-1}{k_2 + 1}, \quad \mbox{and} \quad \gamma = -\alpha \beta.
$$
In particular, we are interested in the case when $k_1 \to 0, k_2 \to + \infty$ (or $\alpha \to -1, \beta \to 1$). Therefore, we may restrict $\gamma > \frac{1}{2}$. First we assume that $r_1 = r_2 = 1$.

\begin{theorem}\label{main_thm}
Let $\varepsilon\in (0,1/2)$ and $\mu \in (0,1)$ be constants. Assume that $u$ is a weak solution of \eqref{main} in $B_1: = B_1(0)$ with $r_1 = r_2 =1$, $k_1 \in (0,1)$, $k_2 \in (1,\infty)$, and  $\gamma \in (1/2, 1)$. For any $m \in \bN$, if $f$ is piecewise $C^{2m-1, \mu}$ in $B_1$, and for some constant $C_m > 0$,
$$
\|u\|_{L^2(B_1)} \le C_m, \quad \|f\|_{C^{2m-1,\mu}(B_1 \cap \cB_j)} \le C_m k_j, \quad j=0,1,2,
$$
then we have
\begin{equation}\label{main_estimate}
|D^m u(x)| \le CC_m \quad \mbox{in}~~B_{1/2},
\end{equation}
where $C > 0$ is a positive constant depending only on $m$ and $\mu$, and in particular is independent of $\varepsilon$, $k_1$, and $k_2$.
\end{theorem}

When $u$ satisfies \eqref{main} in a domain that contains $\cB_1$ and $\cB_2$, we have more precise estimates as follows.

\begin{theorem}\label{global_thm}
Let $\varepsilon\in (0,1/2)$ and $\mu \in (0,1)$ be constants. Assume that $\cB_1 \cup \cB_2 \Subset \cD_1 \Subset \cD$ for some domain $\cD_1$, $u$ is a weak solution of \eqref{main} in $\cD$ with $r_1 = r_2 = 1$, $0 < k_1 < 1$, $k_2 > 1$, and  $\gamma \in (1/2, 1)$. For any $m \in \bN$, if $f$ is piecewise $C^{2m-1, \mu}$ in $\cD$, and for some constant $C_m > 0$,
$$
\|u\|_{L^2(\cD)} \le C_m, \quad \|f\|_{C^{2m-1,\mu}(B_1 \cap \cB_j)} \le C_m \min\{1,k_j\}, \quad j=0,1,2,
$$
then we have
\begin{equation*}
|D^m u(x)| \le \left\{
\begin{aligned}
&CC_m  &&\mbox{in}~~\cD_1 \cap \cB_0,\\
&\frac{CC_m}{k_1+1} &&\mbox{in}~~\cB_1,\\
&\frac{CC_m}{k_2+1} &&\mbox{in}~~\cB_2,
\end{aligned}
\right.
\end{equation*}
where $C>0$ is a constant depending only on  $m$, $\mu$, $\cD_1$, and $\cD$.
\end{theorem}

For the general case when $r_1$ and $r_2$ are not necessarily equal to $1$, we have the following theorem.

\begin{theorem}\label{general_thm}
Let $\varepsilon\in (0,1/2)$, $\mu \in (0,1)$ be constants, and $1/2 < r_1, r_2 < 10$. Then there exist domains $\cD_1 \Subset \cD$ that depend on $r_1, r_2$, such that if $\cB_1 \cup \cB_2 \Subset \cD_1$, and $u$ is a weak solution of \eqref{main} in $\cD$ with $0 < k_1 < 1$, $k_2 > 1$, and  $\gamma \in (1/2, 1)$, we have, for any $m \in \bN$,  if $f$ is piecewise $C^{2m-1, \mu}$ in $\cD$, and for some constant $C_m > 0$,
$$
\|u\|_{L^2(\cD)} \le C_m, \quad \|f\|_{C^{2m-1,\mu}(B_1 \cap \cB_j)} \le C_m \min\{1,k_j\}, \quad j=0,1,2,
$$
then
\begin{equation*}
|D^m u(x)| \le \left\{
\begin{aligned}
&CC_m  &&\mbox{in}~~\cD_1 \cap \cB_0,\\
&\frac{CC_m}{k_1+1} &&\mbox{in}~~\cB_1,\\
&\frac{CC_m}{k_2+1} &&\mbox{in}~~\cB_2,
\end{aligned}
\right.
\end{equation*}
where $C>0$ is a constant depending only on $m$, $\mu$, $r_1$, and $r_2$.
\end{theorem}

In the final part of this paper, we partially answer a question raised by Kang in \cite{Kang}, where in the conclusion section, he mentioned that the extensions to general shape of inclusions and higher dimensions for the case $(k_1-1)(k_2-1)<0$ are quite challenging. We prove the derivatives estimates for the extreme case when $k_1 = 0$, $k_2 = \infty$, and $f = 0$ for general strictly convex inclusions in dimensions $n \ge 2$.

The setting of this problem is as follows. Let $\cD \subset \bR^n$ be a bounded domain containing two smooth relatively strictly convex open sets $\cB_{1}$ and $\cB_{2}$ so that $\cB_1 \cup \cB_2 \Subset \cD$, $\dist(\cB_1 \cap \cB_2) = \varepsilon$, and $\dist(\cB_1 \cup \cB_2, \cD) > 0$. It is known that when $k_1 = 0$, $k_2 = \infty$, and $f = 0$, \eqref{main} is reduced to
\begin{equation*}
\begin{cases}
\Delta{u}=0&\mbox{in}~\cD\setminus\overline{(\cB_1 \cup \cB_2)},\\
\frac{\partial u}{\partial \nu} = 0 &\mbox{on}~\partial{\cB}_{1},\\
u=C \mbox{ (Constant)}&\mbox{on}~\partial{\cB}_{2},\\
\int_{\partial{\cB}_{2}}\frac{\partial{u}}{\partial\nu}=0.\\
\end{cases}
\end{equation*}
We use the notation $x = (x', x_n)$, where $x' \in \bR^{n-1}$. After choosing a coordinate system properly, we can assume that near the origin,  the part of $\partial \cB_1$ and $\partial \cB_2$, denoted by $\Gamma_+$ and $\Gamma_-$, are respectively the graphs of two smooth functions in terms of $x'$. That is,
\begin{align*}
\Gamma_+ = \left\{ x_n = \frac{\varepsilon}{2}+h_1(x'),~|x'|<1\right\}~~ \mbox{and} ~~\Gamma_- = \left\{ x_n = -\frac{\varepsilon}{2}+h_2(x'),~|x'|<1\right\},
\end{align*}
where $h_1$ and $h_2$ are sufficiently smooth functions satisfying
\begin{equation*}
h_1(x')>h_2(x')\quad\mbox{for}~~0<|x'|<1,
\end{equation*}
\begin{equation*}
h_1(0')=h_2(0')=0,\quad D_{x'}h_1(0')=D_{x'}h_2(0')=0, \quad D^2 (h_1-h_2)(0') > 0.
\end{equation*}
For $0 < r\leq 1$, we denote
\begin{align*}
\Omega_{r}:=\left\{(x',x_{n})\in \cD\setminus\overline{(\cB_1 \cup \cB_2)}~\big|~-\frac{\varepsilon}{2}+h_2(x')<x_{n}<\frac{\varepsilon}{2}+h_1(x'),~|x'|<r\right\}.
\end{align*}
We focus on the following problem near the origin:
\begin{equation}\label{main_problem_narrow}
\left\{
\begin{aligned}
-D_i (a^{ij}(x)  D_j u(x)) =0 \quad &\mbox{in }\Omega_{1},\\
a^{ij}(x)  D_j u(x) \nu_i = 0 \quad &\mbox{on } \Gamma_+,\\
u=C \mbox{ (Constant)} \quad &\mbox{on}~\Gamma_-,\\
\end{aligned}
\right.
\end{equation}
where $(a^{ij}(x))$ satisfies, for some constants $\sigma\in (0,1)$ and any $x\in \Omega_1, \xi\in \bR^n$,
\begin{equation*}
\sigma |\xi|^2 \le a^{ij}(x)\xi_i\xi_j, \quad |a^{ij}(x)|\le  \frac{1}{\sigma}.
\end{equation*}

\begin{theorem}\label{extreme_thm}
Assume the above and let $u \in H^1(\Omega_1)$ be a weak solution of \eqref{main_problem_narrow}. For $m \in \bN$, some constants $\alpha \in (0,1)$, and $C_{m,\alpha} > 0$, if
$$
\|a\|_{C^{m-1,\alpha}(\Omega_1)} + \|h_1\|_{C^{m,\alpha}(\{|x'|<1\})} + \|h_2\|_{C^{m,\alpha}(\{|x'|<1\})} \le C_{m,\alpha},
$$
then there exist constants $\mu \in (0,1)$ and $C$, depending only on $n$, $\sigma$, $m$, $\alpha$, and $C_{m,\alpha}$ such that
$$
|D^m u(x)| \le C \mu^{\frac{1}{\sqrt{\varepsilon} + |x'|}} \|u\|_{L^2(\Omega_1)} \quad \mbox{for}~~x \in \Omega_{1/2}.
$$
\end{theorem}

The rest of this paper is organized as follows. First, we review the Green function of the operator $L_{\varepsilon;1,1}$ constructed in \cite{DL}, and derive some preliminary estimates in Section 2. In Section 3, we prove Theorem \ref{main_thm} with $m=1$ to illustrate the main idea without getting into too much technicalities. Then we prove Theorem \ref{main_thm} with general $m \in \bN$ in Section 4, and Theorem \ref{global_thm} in Section 5. In Section 6, we prove Theorem \ref{general_thm} by introducing a conformal map to reduce the problem to the case considered in Theorem \ref{global_thm}. Finally, Theorem \ref{extreme_thm} is proved in Section 7.

\section{Preliminary}
In this section, we first review the Green function of the operator $L_{\varepsilon;1,1}$ constructed in \cite{DL}, and then derive some preliminary estimates. %This section follows closely from Sections 2 and 3 of \cite{DL}.

Let $\Phi_1(x), \Phi_2(x)$ denote the inversion maps of a point $x \in \bR^2$ with respect to $\partial \cB_1$ and $\partial \cB_2$, respectively, that is
$$
\Phi_1(x_1,x_2):= \left( \frac{x_1 - (1+\varepsilon/2)}{(x_1-1-\varepsilon/2)^2 + x_2^2} + 1 + \varepsilon/2, \frac{x_2}{(x_1-1-\varepsilon/2)^2 + x_2^2} \right)
$$
and
$$
\Phi_2(x_1,x_2):= \left( \frac{x_1 +1+\varepsilon/2}{(x_1+1+\varepsilon/2)^2 + x_2^2} - 1 - \varepsilon/2, \frac{x_2}{(x_1+1+\varepsilon/2)^2 + x_2^2} \right).
$$
The auxiliary function $\cG(x,y)$ is given as follows:
\begin{enumerate}[(1)]
\item When $y \in \cB_0$, $\cG(x,y)$ equals
\begin{align*}
&\frac{2}{k_1+1} \sum_{k=0}^\infty (\alpha\beta)^k \Big( \log|(\Phi_1 \Phi_2)^k(x) - y|  - \beta \log|(\Phi_2 \Phi_1)^k \Phi_2(x) - y| \Big)\\
&\hspace{3.5in}\mbox{for}~~x\in\overline{\cB}_1;\\
&\log |x-y| + \sum_{k=0}^\infty \Big[ (\alpha\beta)^k \Big( \log|(\Phi_1 \Phi_2)^k(x) - y| + \log|(\Phi_2 \Phi_1)^k(x) - y| \Big)\\
&-(\alpha\beta)^{k-1} \Big(\beta \log|(\Phi_2 \Phi_1)^k\Phi_2(x) - y|  + \alpha \log|(\Phi_1 \Phi_2)^k\Phi_1(x) - y| \Big)\Big]\\
&\hspace{3.5in}\mbox{for}~~x\in \cB_0;\\
&\frac{2}{k_2+1} \sum_{l=0}^\infty (\alpha\beta)^k \Big( \log|(\Phi_2 \Phi_1)^k(x) - y|  - \alpha \log|(\Phi_1 \Phi_2)^k \Phi_1(x) - y| \Big)\\
&\hspace{3.5in}\mbox{for}~~x\in\overline{\cB}_2;\\
\end{align*}
\item When $y \in \cB_1$, $\cG(x,y)$ equals
\begin{align*}
&\frac{1}{k_1}(\log|x-y| + \alpha \log |\Phi_1(x) - y|) - \frac{4\beta}{(k_1 + 1)^2} \sum_{k=0}^\infty (\alpha\beta)^k \log|(\Phi_2 \Phi_1)^k \Phi_2(x) - y|\\
&\hspace{3in}\mbox{for}~~x \in \overline{\cB}_1 \setminus \{(1+\varepsilon/2,0)\};\\
&\frac{2}{k_1+1} \sum_{k=0}^\infty (\alpha\beta)^k \Big( \log|(\Phi_2 \Phi_1)^k(x) - y|  - \beta \log|(\Phi_2 \Phi_1)^k \Phi_2(x) - y| \Big)\\
&\hspace{3in}\mbox{for}~~x\in\cB_0;\\
&\frac{4}{(k_1+1)(k_2+1)}\sum_{k=0}^\infty (\alpha\beta)^k \log|(\Phi_2 \Phi_1)^k(x) - y|\\
&\hspace{3in}\mbox{for}~~x\in\overline{\cB}_2;\\
\end{align*}
\item When $y \in \cB_2$, $\cG(x,y)$ equals
\begin{align*}
&\frac{4}{(k_1+1)(k_2+1)}\sum_{k=0}^\infty (\alpha\beta)^k \log|(\Phi_1 \Phi_2)^k(x) - y|\\
&\hspace{3in}\mbox{for}~~x\in\overline{\cB}_1;\\
&\frac{2}{k_2+1} \sum_{k=0}^\infty (\alpha\beta)^k \Big( \log|(\Phi_1 \Phi_2)^k(x) - y|  - \alpha \log|(\Phi_1 \Phi_2)^k \Phi_1(x) - y| \Big)\\
&\hspace{3in}\mbox{for}~~x\in\cB_0;\\
&\frac{1}{k_2}(\log|x-y| + \beta \log |\Phi_2(x) - y|) - \frac{4\alpha}{(k_2 + 1)^2} \sum_{k=0}^\infty (\alpha\beta)^k \log|(\Phi_1 \Phi_2)^k \Phi_1(x) - y|\\
&\hspace{3in}\mbox{for}~~x\in\overline{\cB}_2 \setminus \{(-1-\varepsilon/2,0)\}.\\
\end{align*}
\end{enumerate}
Define
\begin{equation}\label{green_function}
G(x,y) = \left\{
\begin{aligned}
&\cG(x,y) &&\mbox{for}~~y \in \cB_0,\\
&\cG(x,y)+\frac{\alpha}{1-\alpha}\cG(x,(1+\varepsilon/2,0)) &&\mbox{for}~~y \in \cB_1,\\
&\cG(x,y)+\frac{\beta}{1-\beta}\cG(x,(-1-\varepsilon/2,0)) &&\mbox{for}~~y \in \cB_2.
\end{aligned}
\right.
\end{equation}
By \cite{DL}*{Proposition 2.3}, $G$ is a Green function of $L_{\varepsilon;1,1}$ in the sense that
$$
a(x) \Delta_x G(x,y) = \delta(x-y) \quad \mbox{for}~~ x \not\in \partial \cB_1 \cup \partial \cB_2,
$$
and $G(\cdot, y)$, $aD_\nu G(\cdot,y)$ are continuous across $\partial\cB_1$ and $\partial\cB_2$. Let $\eta \in C^\infty_0(B_{3/4})$ be a cutoff function such that $\eta = 1$ in $B_{1/2}$. Let $v = u \eta$, where $u$ is a solution of \eqref{main}. Then $v$ satisfies
$$
D_i(a(x) D_i v) = D_i \tilde f_i + \tilde{f}_3 \quad \mbox{in}~~\bR^2,
$$
where
\begin{equation}\label{tilde_f_formula}
\tilde f_i = f_i \eta + au D_i \eta, \quad \tilde{f}_3 = - f_i D_i \eta + a D_i u D_i \eta.
\end{equation}
We define
\begin{align}
\tilde{u}(x) = &- \int_{\cB_1} D_{y_i} G(x,y) \tilde{f}_i(y) \, dy - \int_{\cB_2} D_{y_i} G(x,y) \tilde{f}_i(y) \, dy\nonumber\\
& - \int_{\cB_0} D_{y_i} G(x,y) \tilde{f}_i(y) \, dy + \int_{B_1}  G(x,y) \tilde f_3(y) \, dy\nonumber\\
:= & - w_1(x) - w_2(x) - w_0(x) + w_3(x). \label{tilde_u_def}
\end{align}
We know from \cite{DL}*{pp. 1447} that $u = \tilde{u} + C_0$ for some constant $C_0$.
We define for $j = 0,1,2,$
\begin{equation}\label{h_def}
h_j(x) = \int_{\cB_j} D_{y_i} \log |x-y| \tilde{f}_i(y) \, dy
\end{equation}
and
\begin{equation}\label{g_def}
g_j(x) = \int_{\cB_j} \log |x-y| \tilde{f}_3(y) \, dy.
\end{equation}
Since supp$(D \eta) \subset B_{3/4}\setminus B_{1/2}$, by \cite{DL}*{Lemma 3.2} and \cite{DZ}*{Lemma 2.1}, we have for $m \in \bN$,
\begin{equation}
\| u \|_{C^{2m,\mu}(\cB_j\cap \supp(D \eta))}  \le CC_m, \,\, \| \tilde{f}_i \|_{C^{2m-1,\mu}(\cB_j)}\le CC_m k_j \,\, \mbox{for}~~i=1,2,3,~~ j=0,1,2,\label{uf_derivatives}
\end{equation}
and
\begin{align}
\|h_j\|_{C^{2m,\mu}(B_3 \cap \cB_i)} + \|g_j\|_{C^{2m,\mu}(B_3 \cap \cB_i)}\le CC_m k_j \quad \mbox{for}~~i,j = 0,1,2 .\label{hg_derivatives}
\end{align}

From now on, we identify a point $x = (x_1,x_2) \in \bR^2$ with a complex number $z = x_1 + ix_2 \in \bC$. We will derive some derivative estimates of the maps $(\Phi_2 \Phi_1)^k(z)$ with respect to $k$ and $z \in \cB_2$. Estimates of $(\Phi_1 \Phi_2)^k(z)$ with respect to $k$ and $z$ for $z \in \cB_1$ will follow similarly.

We denote $a = 1 + \varepsilon/2$ for convenience, after a change of variable
$$
2az - (2a^2-1) \rightarrow z,
$$
we have
$$
(\Phi_2 \Phi_1)(z) = -1/z - 2(2a^2 -1),
$$
and the two fixed points of $\Phi_2\Phi_1$ are given by
\begin{align*}
\lambda_1:&= - (2a^2 - 1) + 2a \sqrt{a^2-1} \sim -1 + 2 \sqrt{\varepsilon},\\
\lambda_2:&= - (2a^2 - 1) - 2a \sqrt{a^2-1} \sim -1 - 2 \sqrt{\varepsilon}.
\end{align*}
We denote $ \psi = \Phi_2 \Phi_1$ for simplicity. For $z \in \cB_2$, we denote $r:= |z - \lambda_2|,$ and
$$
I_{k}:= (z - \lambda_2^{-1}) - |\lambda_2|^{-k}(z-\lambda_2) = (z - \lambda_2^{-1})(1-|\lambda_2|^{-k}) + (\lambda_2 - \lambda_2^{-1})|\lambda_2|^{-k}.
$$
It is easy to see that $|z - \lambda_2^{-1}| \sim r + \sqrt{\varepsilon}$, Re $(z - \lambda_2^{-1}) < 0$, and $\lambda_2 - \lambda_2^{-1} \sim - \sqrt{\varepsilon}$. Hence
$$
|I_{k}| \sim (r + \sqrt{\varepsilon})(1 - \lambda_2^{-k}) + \sqrt{\varepsilon}.
$$
Following the iteration argument from Section 3 of \cite{DL}, we have, for any $\alpha \in \bN$,
\begin{align}
\psi^k (z)&= \lambda_2 + (\lambda_2^2 -1) \lambda_2^{-2k-1} (z - \lambda_2) I_{2k}^{-1},\label{psi_k}\\
D^\alpha (\psi^k (z)) &= \frac{(\lambda_2 - \lambda_2^{-1})^2}{\lambda_2^{2k}} (-1)^{\alpha - 1}\alpha! (1 - \lambda_2^{-2k})^{\alpha - 1}I_{2k}^{-(\alpha+1)},\label{dx}
\end{align}
and in particular, since $|\lambda_2 - \lambda_2^{-1}| \lesssim \sqrt{\varepsilon}$ and $|I_2| \gtrsim \sqrt{\varepsilon}$,
\begin{equation}\label{D_psi}
|D^\alpha (\psi (z))| \le C,
\end{equation}
where $C$ is a positive constant depending only on $\alpha$.
By \eqref{psi_k}, for any $\beta \ge 1$,
\begin{align*}
D_k^\beta(\psi^{k}) (z) &= D_k^\beta \Big[ (\lambda_2^2 -1) \lambda_2^{-2k-1} (z - \lambda_2) I_{2k}^{-1} \Big]  \\
&=  \sum_{\beta_1 + \beta_2 = \beta} C_{\beta_1, \beta_2} (\lambda_2^2 -1)(z - \lambda_2) [D_k^{\beta_1} \lambda_2^{-2k-1}] [D_k^{\beta_2 }I_{2k}^{-1}].
\end{align*}
Since
$$
D_k \lambda_2^{-2k} = -2 \log |\lambda_2| \lambda_2^{-2k} \sim - \sqrt{\varepsilon} \lambda_2^{-2k}
$$
and
$$
D_k^{\beta_2 }I_{2k}^{-1} = \sum_{j=0}^{\beta_2} C_j  (\log |\lambda_2|)^{\beta_2} \lambda_2^{-2jk} (z - \lambda_2)^j I_{2k}^{-j-1},
$$
where $C_j$ is some constant independent of $k$ and $\varepsilon$, we obtain
\begin{equation}\label{dk}
|D_k^\beta(\psi^{k}) (z)| \lesssim \frac{\varepsilon^{(1 + \beta)/2}}{\lambda_2^{2k}} \sum_{j = 0}^\beta \frac{r^{j+1}}{|I_{2k}|^{j+1}}.
\end{equation}
For any $\beta \ge 0$ and $\alpha \in \bN$, by \eqref{dx}, we have
\begin{align*}
D^\beta_k D^\alpha (\psi^k (z)) = \sum_{\beta_1 + \beta_2 + \beta_3 = \beta} & C_{\beta_1,\beta_2,\beta_3} (\lambda_2 - \lambda_2^{-1})^2 (-1)^{\alpha - 1}\alpha! [D_k^{\beta_1} \lambda_2^{-2k}] \times \\
&\times [D_k^{\beta_2} (1 - \lambda_2^{-2k})^{\alpha - 1}] [D_k^{\beta_3} I_{2k}^{-(\alpha+1)}].
\end{align*}
Observe that
\begin{align*}
\left| D_k^{\beta_2} (1 - \lambda_2^{-2k})^{\alpha - 1} \right| &= \left| \sum_{j=0}^{\min(\beta_2 , \alpha-1)} C_j (1 - \lambda_2^{-2k})^{\alpha -j - 1} \lambda_2^{-2jk} (\log |\lambda_2|)^{\beta_2} \right|\\
& \lesssim \varepsilon^{\beta_2 /2} \sum_{j=0}^{\min(\beta_2 , \alpha-1)} |1 - \lambda_2^{-2k}|^{\alpha -j - 1},
\end{align*}
and
\begin{align*}
|D_k^{\beta_3} I_{2k}^{-(\alpha+1)}| &= \left|\sum_{j=0}^{\beta_3} C_j (\log |\lambda_2|)^{\beta_3} \lambda_2^{-2jk}(z - \lambda_2)^{j}  I_{2k}^{-(\alpha+j+1)}\right|\\
& \lesssim \varepsilon^{\beta_3/2} \sum_{j=0}^{\beta_3} \frac{r^j}{|I_{2k}|^{\alpha+j+1}},
\end{align*}
where $C_j$ is some constant and can be different from line to line. Thus we have
\begin{align}\label{dkdx}
|D^\beta_k D^\alpha (\psi^k (z))| \lesssim \frac{\varepsilon^{1 + \beta/2}}{\lambda_2^{2k}} \sum_{\substack{ \beta_2 + \beta_3 \le \beta\\
\beta_2 \le \alpha -1}} |1 - \lambda_2^{-2k}|^{\alpha - \beta_2 - 1} \frac{r^{\beta_3}}{|I_{2k}|^{\alpha + \beta_3 + 1}}.
\end{align}
Note that $\beta_2, \beta_3$ in \eqref{dkdx} might be different from the ones above.

\section{First order derivative}
In this section, we prove Theorem \ref{main_thm} when $m = 1$. Without loss of generality, we may assume that $C_1 = 1$.

\textbf{Case 1:} $x \in \cB_0 \cap B_{1/2}$. In this case, $\psi^k(x) \in \cB_2$ for $k \in \bN$. By \eqref{green_function} and \eqref{tilde_u_def}, we have
\begin{align*}
Dw_1(x) = &\frac{2}{k_1+1} \sum_{k=0}^\infty (-1)^k \gamma^k \Big( D[h_1(\psi^k (x))] - \beta D[h_1(\psi^k \Phi_2 (x))] \Big).
\end{align*}
We denote
$$
\Theta_k:= \gamma^k D[h_1(\psi^k (x))]= \gamma^k D h_1 [\psi^{k} (x)]D(\psi^{k-1}) [\psi(x)] D\psi (x).
$$
Then
$$
\sum_{k=0}^\infty (-1)^k \Theta_k = \frac{1}{2} \sum_{k=1}^\infty (-1)^k (\Theta_k - \Theta_{k+1}) + \Theta_0 - \frac{1}{2} \Theta_1.
$$
By \eqref{uf_derivatives} and \eqref{hg_derivatives}, it is clear that $ |\Theta_0 - \frac{1}{2} \Theta_1| \le C k_1$, where $C$ is a positive constant independent of $\varepsilon$ and $\gamma$. Therefore, the goal is to estimate $|\Theta_k - \Theta_{k+1}|$. By the mean value theorem,
$$
\Theta_k - \Theta_{k+1} =  D_k \Theta_{k}\Big|_{k = \bar{k}}
$$
for some $\bar{k} \in (k, k+1)$. We claim that
\begin{equation}\label{goal_1}
 \sum_{k = 1}^\infty |\Theta_k - \Theta_{k+1}| \lesssim \sum_{k = 1}^\infty |D_k \Theta_{k}| \lesssim k_1.
\end{equation}
By the chain rule and the product rule, we have
\begin{align*}
D_k \Theta_k =& (\log \gamma)\gamma^k D h_1 [\psi^{k} (x)]D(\psi^{k-1}) [\psi(x)] D\psi (x)\\
&+ \gamma^k  D^2 h_1 [\psi^{k} (x)] D_k(\psi^{k-1}) [\psi(x)] D(\psi^{k-1}) [\psi(x)] D\psi (x)\\
&+ \gamma^k D h_1 [\psi^{k} (x)]D_kD(\psi^{k-1}) (x)D\psi (x)\\
=&: J_k^1 + J_k^2 + J_k^3.
\end{align*}
Recall that $|\lambda_2| > 1$ and $|\log \gamma| \sim 1-\gamma$, by \eqref{hg_derivatives}, \eqref{D_psi}, \eqref{dk}, and \eqref{dkdx}, we have
\begin{align*}
|J_k^1|\lesssim& k_1 |\log \gamma| \gamma^k \frac{\varepsilon}{\lambda_2^{2k} |I_{2k}|^2},\\
|J_k^2|\lesssim& k_1 \gamma^k \frac{\varepsilon}{\lambda_2^{2k}} \left( \frac{r}{|I_{2k}|} + \frac{r^2}{|I_{2k}|^2} \right) \frac{\varepsilon}{\lambda_2^{2k}} \frac{1}{|I_{2k}|^2}
\lesssim k_1 \frac{\varepsilon^2}{\lambda_2^{2k}} \left(\frac{r}{|I_{2k}|^3} + \frac{r^2}{|I_{2k}|^4} \right),\\
|J_k^3|\lesssim& k_1 \gamma^k \frac{\varepsilon^{3/2}}{\lambda_2^{2k}} \left(\frac{1}{|I_{2k}|^2} + \frac{r}{|I_{2k}|^3} \right).
\end{align*}
When $k < \frac{1}{r+\sqrt{\varepsilon}}$, we have
$$
|I_{2k}| \gtrsim \sqrt{\varepsilon}.
$$
Therefore,
\begin{align*}
\Big| \sum_{k < \frac{1}{r+\sqrt{\varepsilon}}} J_k^1 \Big| & \lesssim k_1 (1 - \gamma)\sum_{k < \frac{1}{r+\sqrt{\varepsilon}}}  \gamma^k \frac{\varepsilon}{\lambda_2^{2k} |I_{2k}|^2}\\
&\lesssim k_1 (1-\gamma) \sum_{k < \frac{1}{r+\sqrt{\varepsilon}}} \gamma^k \lesssim k_1,
\end{align*}
\begin{align*}
\Big| \sum_{k < \frac{1}{r+\sqrt{\varepsilon}}} J_k^2 \Big| &\lesssim k_1 \sum_{k < \frac{1}{r+\sqrt{\varepsilon}}} \frac{\varepsilon^2}{\lambda_2^{2k}} \left(\frac{r}{|I_{2k}|^3} + \frac{r^2}{|I_{2k}|^4} \right)\\
& \lesssim k_1 \sum_{k < \frac{1}{r+\sqrt{\varepsilon}}}\varepsilon^2 \left(\frac{r}{\varepsilon^{3/2}} + \frac{r^2}{\varepsilon^{2}}\right) \\
& \lesssim k_1 \frac{\varepsilon^2}{r + \sqrt{\varepsilon}}\left(\frac{r}{\varepsilon^{3/2}} + \frac{r^2}{\varepsilon^{2}} \right) \lesssim k_1 (\sqrt{\varepsilon} + r),
\end{align*}
and
\begin{align*}
\Big| \sum_{k < \frac{1}{r+\sqrt{\varepsilon}}} J_k^3 \Big| &\lesssim k_1 \sum_{k < \frac{1}{r+\sqrt{\varepsilon}}}\gamma^k \frac{\varepsilon^{3/2}}{\lambda_2^{2k}} \left(\frac{1}{|I_{2k}|^2} + \frac{r}{|I_{2k}|^3} \right)\\
& \lesssim k_1 \sum_{k < \frac{1}{r+\sqrt{\varepsilon}}}\varepsilon^{3/2} \left(\frac{1}{\varepsilon} + \frac{r}{\varepsilon^{3/2}} \right)\\
& \lesssim k_1 \frac{\varepsilon^{3/2}}{r + \sqrt{\varepsilon}}\left(\frac{1}{\varepsilon} + \frac{r}{\varepsilon^{3/2}} \right) \lesssim k_1.
\end{align*}
When $\frac{1}{r+ \sqrt{\varepsilon}} \le k < \frac{1}{\sqrt{\varepsilon}}$, we have
$$
|I_{2k}| \gtrsim (r+\sqrt{\varepsilon})k\sqrt{\varepsilon}.
$$
Therefore,
\begin{align*}
\Big| \sum_{\frac{1}{r+ \sqrt{\varepsilon}} \le k < \frac{1}{\sqrt{\varepsilon}}} J_k^1 \Big| & \lesssim k_1 \sum_{\frac{1}{r+ \sqrt{\varepsilon}} \le k < \frac{1}{\sqrt{\varepsilon}}} |\log \gamma| \gamma^k \frac{\varepsilon}{\lambda_2^{2k} |I_{2k}|^2}\\
&\lesssim k_1 \sum_{\frac{1}{r+ \sqrt{\varepsilon}} \le k < \frac{1}{\sqrt{\varepsilon}}} |\log \gamma| \gamma^k \frac{\varepsilon}{\lambda_2^{2k} (r+\sqrt{\varepsilon})^2 k^2 \varepsilon}\\
&\lesssim k_1 \gamma^{\frac{1}{r+\sqrt{\varepsilon}}} |\log \gamma| \frac{1}{(r+ \sqrt{\varepsilon})^2} \sum_{\frac{1}{r+ \sqrt{\varepsilon}} \le k < \frac{1}{\sqrt{\varepsilon}}} \frac{1}{k^2}\\
&\lesssim k_1 e^{\frac{\log \gamma}{r+ \sqrt{\varepsilon}}}  \frac{|\log \gamma|}{r+ \sqrt{\varepsilon}} \lesssim k_1,
\end{align*}
where we used the fact that $|e^{-x}x| \le C$ for any $x > 0$,
\begin{align*}
\Big| \sum_{\frac{1}{r+ \sqrt{\varepsilon}} \le k < \frac{1}{\sqrt{\varepsilon}}} J_k^2 \Big| & \lesssim k_1 \sum_{\frac{1}{r+ \sqrt{\varepsilon}} \le k < \frac{1}{\sqrt{\varepsilon}}}\gamma^k \frac{\varepsilon^2}{\lambda_2^{2k}} \left(\frac{r}{|I_{2k}|^3} + \frac{r^2}{|I_{2k}|^4} \right)\\
& \lesssim k_1 \sum_{\frac{1}{r+ \sqrt{\varepsilon}} \le k < \frac{1}{\sqrt{\varepsilon}}} \varepsilon^2 \left( \frac{r}{(r+\sqrt{\varepsilon})^3k^3\varepsilon^{3/2}} + \frac{r^2}{(r+\sqrt{\varepsilon})^4k^4\varepsilon^{2}}\right)\\
& \lesssim k_1 \sum_{\frac{1}{r+ \sqrt{\varepsilon}} \le k < \frac{1}{\sqrt{\varepsilon}}} \left( \frac{1}{(r+\sqrt{\varepsilon})k^3} + \frac{1}{(r+\sqrt{\varepsilon})^2k^4} \right)\\
& \lesssim k_1( r + \sqrt{\varepsilon} ),
\end{align*}
and
\begin{align*}
\Big| \sum_{\frac{1}{r+ \sqrt{\varepsilon}} \le k < \frac{1}{\sqrt{\varepsilon}}} J_k^3 \Big| & \lesssim k_1 \sum_{\frac{1}{r+ \sqrt{\varepsilon}} \le k < \frac{1}{\sqrt{\varepsilon}}} \gamma^k \frac{\varepsilon^{3/2}}{\lambda_2^{2k}} \left(\frac{1}{|I_{2k}|^2} + \frac{r}{|I_{2k}|^3} \right)\\
& \lesssim k_1 \sum_{\frac{1}{r+ \sqrt{\varepsilon}} \le k < \frac{1}{\sqrt{\varepsilon}}} \varepsilon^{3/2} \left( \frac{1}{(r+\sqrt{\varepsilon})^2k^2\varepsilon} + \frac{r}{(r+\sqrt{\varepsilon})^3k^3\varepsilon^{3/2}}\right)\\
& \lesssim k_1 \sum_{\frac{1}{r+ \sqrt{\varepsilon}} \le k < \frac{1}{\sqrt{\varepsilon}}} \left( \frac{1}{(r+\sqrt{\varepsilon})k^2} + \frac{1}{(r+\sqrt{\varepsilon})^2k^3} \right) \lesssim k_1.
\end{align*}
Finally, when $k \ge \frac{1}{\sqrt{\varepsilon}}$, we have
$$
|I_{2k}| \gtrsim r+\sqrt{\varepsilon}.
$$
Therefore,
\begin{align*}
\Big| \sum_{k \ge \frac{1}{\sqrt{\varepsilon}}} J_k^1 \Big| & \lesssim k_1 \sum_{k \ge \frac{1}{\sqrt{\varepsilon}}} |\log \gamma| \gamma^k \frac{\varepsilon}{\lambda_2^{2k} |I_{2k}|^2}\\
&\lesssim  k_1(1-\gamma) \sum_{k \ge \frac{1}{\sqrt{\varepsilon}}} \gamma^k \lesssim k_1 \gamma^{\frac{1}{\sqrt{\varepsilon}}},
\end{align*}
\begin{align*}
\Big| \sum_{k \ge \frac{1}{\sqrt{\varepsilon}}} J_k^2 \Big| & \lesssim k_1 \sum_{k \ge \frac{1}{\sqrt{\varepsilon}}} \gamma^k \frac{\varepsilon^2}{\lambda_2^{2k}} \left(\frac{r}{|I_{2k}|^3} + \frac{r^2}{|I_{2k}|^4} \right)\\
&\lesssim k_1 \sum_{k \ge \frac{1}{\sqrt{\varepsilon}}} \frac{\gamma^k}{\lambda_{2}^{2k}} \varepsilon^2  \left( \frac{r}{(r+\sqrt{\varepsilon})^3} + \frac{r^2}{(r+\sqrt{\varepsilon})^4}\right)\\
&\lesssim k_1 \frac{\varepsilon (\gamma/\lambda_2^2)^{\frac{1}{\sqrt{\varepsilon}}}}{1- \gamma/\lambda_2^2} \lesssim k_1 \sqrt{\varepsilon} (\gamma/\lambda_2^2)^{\frac{1}{\sqrt{\varepsilon}}},
\end{align*}
where we used $1- \gamma/\lambda_2^2 \gtrsim \sqrt{\varepsilon}$, and
\begin{align*}
\Big| \sum_{k \ge \frac{1}{\sqrt{\varepsilon}}} J_k^3 \Big| & \lesssim k_1 \sum_{k \ge \frac{1}{\sqrt{\varepsilon}}} \gamma^k \frac{\varepsilon^{3/2}}{\lambda_2^{2k}} \left(\frac{1}{|I_{2k}|^2} + \frac{r}{|I_{2k}|^3} \right)\\
&\lesssim k_1 \sum_{k \ge \frac{1}{\sqrt{\varepsilon}}} \frac{\gamma^k}{\lambda_{2}^{2k}} \varepsilon^{3/2}  \left( \frac{1}{(r+\sqrt{\varepsilon})^2} + \frac{r}{(r+\sqrt{\varepsilon})^3}\right)\\
&\lesssim  k_1\frac{\sqrt\varepsilon (\gamma/\lambda_2^2)^{\frac{1}{\sqrt{\varepsilon}}}}{1- \gamma/\lambda_2^2} \lesssim k_1 (\gamma/\lambda_2^2)^{\frac{1}{\sqrt{\varepsilon}}}.
\end{align*}
Therefore, \eqref{goal_1} follows from the estimates above, and hence
$$
\left| \sum_{k=0}^\infty (-1)^k \Theta_k \right| \le Ck_1,
$$
where $C$ is a positive constant independent of $\varepsilon$ and $\gamma$.

By the same argument, we can estimate
$$
\left| \sum_{k=0}^\infty (-1)^k \gamma^k D[h_1(\psi^k \Phi_2 (x))] \right| \le Ck_1.
$$
Therefore,
$$
|Dw_1(x)| \le \frac{Ck_1}{k_1 + 1} \quad \mbox{for}~~x \in \cB_0 \cap B_{1/2}.
$$

By \eqref{green_function} and \eqref{tilde_u_def} again, we have, for $x \in \cB_0 \cap B_{1/2}$,
\begin{align*}
D w_2(x) = &\frac{2}{k_2+1} \sum_{k=0}^\infty (-1)^k \gamma^k \Big( D[h_2((\Phi_1\Phi_2)^k (x))] - \alpha D[h_2((\Phi_1\Phi_2)^k \Phi_1 (x))] \Big),\\
D w_0(x) = & D h_0(x) + \sum_{k=0}^\infty \Big[(-1)^k \gamma^k \Big( D[h_0((\Phi_1\Phi_2)^k (x))] + D[h_0((\Phi_2\Phi_1)^k (x))] \Big)\\
& - (-1)^{k-1} \gamma^{k-1} \Big( \beta D [h_0((\Phi_2\Phi_1)^{k-1} \Phi_2 (x))] + \alpha D[h_0((\Phi_1\Phi_2)^{k-1} \Phi_1 (x))] \Big) \Big].
\end{align*}
Therefore, in the same way we can estimate
$$
|Dw_2(x)| \le \frac{Ck_2}{k_2+1} \quad \mbox{and} \quad |Dw_0(x)| \le C \quad \mbox{for}~~x \in \cB_0 \cap B_{1/2},
$$
where $C$ is a positive constant independent of $\varepsilon$, $k_1$, and $k_2$. To estimate $w_3$, since supp$(\tilde{f}_3) \subset B_1$, we can write
\begin{align*}
w_3(x) &= \int_{\cB_1}  G(x,y) \tilde f_3(y) \, dy + \int_{\cB_2}  G(x,y) \tilde f_3(y) \, dy + \int_{\cB_0}  G(x,y) \tilde f_3(y) \, dy\\
&= \frac{2}{k_1+1} \sum_{k=0}^\infty (-1)^k \gamma^k \Big( g_1((\Phi_2\Phi_1)^k (x)) - \beta g_1((\Phi_2\Phi_1)^k \Phi_2 (x)) \Big)\\
&\quad+ \alpha \sum_{k=0}^\infty (-1)^k \gamma^k \Big( \log|(\Phi_2\Phi_1)^k (x) - (a,0)| \\
&\quad- \beta \log|(\Phi_2\Phi_1)^k \Phi_2(x) - (a,0)| \Big) \int_{\cB_1} \tilde f_3(y) \, dy\\
&\quad+ \frac{2}{k_2+1} \sum_{k=0}^\infty (-1)^k \gamma^k \Big( g_2((\Phi_1\Phi_2)^k (x)) - \alpha g_2((\Phi_1\Phi_2)^k \Phi_1 (x)) \Big)\\
&\quad+ \beta \sum_{k=0}^\infty (-1)^k \gamma^k \Big( \log|(\Phi_1\Phi_2)^k (x) + (a,0)| \\
&\quad- \alpha \log|(\Phi_1\Phi_2)^k \Phi_1(x) + (a,0)| \Big) \int_{\cB_2} \tilde f_3(y) \, dy\\
&\quad+ g_0(x) +  \sum_{k=0}^\infty \Big[(-1)^k \gamma^k \Big( g_0((\Phi_1\Phi_2)^k (x)) + g_0((\Phi_2\Phi_1)^k (x)) \Big)\\
&\quad- (-1)^{k-1} \gamma^{k-1} \Big( \beta g_0((\Phi_2\Phi_1)^{k-1} \Phi_2 (x)) + \alpha g_0((\Phi_1\Phi_2)^{k-1} \Phi_1 (x)) \Big) \Big].
\end{align*}
Note that for any $x \in \cB_0 \cap B_{1/2}$,
\begin{align*}
&|(\Phi_2\Phi_1)^k (x) - (a,0)| \ge 1, &&|(\Phi_2\Phi_1)^k \Phi_2(x) - (a,0)| \ge 1,\\
&|(\Phi_1\Phi_2)^k (x) + (a,0)| \ge 1, &&|(\Phi_1\Phi_2)^k \Phi_1(x) + (a,0)| \ge 1.
\end{align*}
By \eqref{tilde_f_formula} and \eqref{main},
$$
\int_{B_1} \tilde{f}_3 = \int_{B_1} -f_i D_i \eta + a(x) D_i u D_i \eta = 0.
$$
Therefore, by \eqref{uf_derivatives},
$$
\left| \int_{\cB_2} \tilde{f}_3 \right| = \left|  \int_{\cB_1 \cup \cB_0} -f_i D_i \eta + a(x) D_i u D_i\eta \right|  \le CC_m,
$$
where $C$ is some positive constant independent of $\varepsilon$, $k_1$, and $k_2$. Then we can estimate $Dw_3$ similarly as above to obtain
$$
|Dw_3(x)| \le C \quad \mbox{for}~~x \in \cB_0 \cap B_{1/2}.
$$
Therefore,
$$
|Du(x)| \le C \quad \mbox{for}~~x \in \cB_0 \cap B_{1/2},
$$
where $C$ is a positive constant independent of $\varepsilon$, $k_1$, and $k_2$.

\textbf{Case 2:}  $x \in \cB_1 \cap B_{1/2}$. By \eqref{green_function} and \eqref{tilde_u_def}, we have
\begin{align*}
 w_1(x) = & \frac{1}{k_1} h_1(x) + \frac{\alpha}{k_1}h_1(\Phi_1(x))- \frac{4\beta}{(k_1+1)^2} \sum_{k=0}^\infty (-1)^k \gamma^k  h_1((\Phi_2\Phi_1)^k \Phi_2 (x)),\\
 w_2(x) = & \frac{4}{(k_1+1)(k_2+1)}\sum_{k=0}^\infty (-1)^k \gamma^k h_2((\Phi_1\Phi_2)^k (x)),\\
 w_0(x) = & \frac{2}{k_1+1} \sum_{k=0}^\infty (-1)^k \gamma^k \Big( h_0((\Phi_1\Phi_2)^k (x)) - \beta h_0((\Phi_2\Phi_1)^k \Phi_2 (x)) \Big),
\end{align*}
and
\begin{align*}
w_3(x) &= \int_{\cB_1}  G(x,y) \tilde f_3(y) \, dy + \int_{\cB_2}  G(x,y) \tilde f_3(y) \, dy + \int_{\cB_0}  G(x,y) \tilde f_3(y) \, dy\\
&=  \frac{1}{k_1}g_1(x) + \frac{\alpha}{k_1}g_1(\Phi_1(x)) - \frac{4\beta}{(k_1+1)^2} \sum_{k=0}^\infty (-1)^k \gamma^k g_1((\Phi_2\Phi_1)^k \Phi_2 (x))\\
&\quad+ \frac{\alpha}{k_1}\log|x-(a,0)| \int_{\cB_1} \tilde f_3(y) \, dy\\
&\quad- \frac{2 \alpha\beta}{k_1+1} \sum_{k=0}^\infty (-1)^k \gamma^k \log |(\Phi_2\Phi_1)^k \Phi_2(x) - (a,0)|\int_{\cB_1} \tilde f_3(y) \, dy\\
&\quad+ \frac{4}{(k_1+1)(k_2+1)}\sum_{k=0}^\infty (-1)^k \gamma^k g_2((\Phi_1\Phi_2)^k (x))\\
&\quad+ \frac{2 \beta}{k_1+1} \sum_{k=0}^\infty (-1)^k \gamma^k \log|(\Phi_1\Phi_2)^k (x)+(a,0)| \int_{\cB_2} \tilde f_3(y) \, dy\\
&\quad+ \frac{2}{k_1+1} \sum_{k=0}^\infty (-1)^k \gamma^k \Big( g_0((\Phi_1\Phi_2)^k (x)) - \beta g_0((\Phi_2\Phi_1)^k \Phi_2 (x)) \Big),
\end{align*}
where we used $\log |x-(a,0)| = - \log |\Phi_1(x) - (a,0)|$ for $x \in \cB_1 \cap B_{1/2}$.
Note that
\begin{align*}
|x-(a,0)| \ge 1/2, \quad |(\Phi_2\Phi_1)^k \Phi_2(x) - (a,0)| \ge 1, \quad |(\Phi_1\Phi_2)^k (x)+(a,0)| \ge 1.
\end{align*}
We can estimate $Dw_i,i=0,1,2,3,$ as above to obtain
\begin{align*}
|Dw_1(x)|&\le C, \quad |Dw_2(x)| \le \frac{Ck_2}{(k_1+1)(k_2+1)},\\
|Dw_0(x)|&\le \frac{C}{k_1+1},\quad|Dw_3(x)| \le C.
\end{align*}
Therefore,
$$
|Du(x)| \le C \quad \mbox{for}~~x \in \cB_1 \cap B_{1/2},
$$
where $C$ is a positive constant independent of $\varepsilon$, $k_1$, and $k_2$.

\textbf{Case 3:} $x \in \cB_2 \cap B_{1/2}$. By \eqref{green_function} and \eqref{tilde_u_def} again, we have
\begin{align*}
w_1(x) = & \frac{4}{(k_1+1)(k_2+1)}\sum_{k=0}^\infty (-1)^k \gamma^k h_1((\Phi_2\Phi_1)^k (x)),\\
w_2(x) = &\frac{1}{k_2} h_2(x) + \frac{\beta}{k_2}h_2(\Phi_2(x))- \frac{4\alpha}{(k_2+1)^2} \sum_{k=0}^\infty (-1)^k \gamma^k  h_2((\Phi_1\Phi_2)^k \Phi_1 (x)),\\
w_0(x) = & \frac{2}{k_2+1} \sum_{k=0}^\infty (-1)^k \gamma^k \Big( h_0((\Phi_2\Phi_1)^k (x)) - \alpha h_0((\Phi_1\Phi_2)^k \Phi_1 (x)) \Big),
\end{align*}
and
\begin{align*}
w_3(x) &= \int_{\cB_1}  G(x,y) \tilde f_3(y) \, dy + \int_{\cB_2}  G(x,y) \tilde f_3(y) \, dy + \int_{\cB_0}  G(x,y) \tilde f_3(y) \, dy\\
&=\frac{4}{(k_1+1)(k_2+1)}\sum_{k=0}^\infty (-1)^k \gamma^k g_1((\Phi_2\Phi_1)^k (x))\\
&\quad + \frac{2\alpha}{k_2+1} \sum_{k=0}^\infty (-1)^k \gamma^k \log|(\Phi_2\Phi_1)^k (x)-(a,0)| \int_{\cB_1} \tilde f_3(y) \, dy\\
&\quad+ \frac{1}{k_2}g_2(x) + \frac{\beta}{k_2}g_2(\Phi_2(x)) - \frac{4\alpha}{(k_2+1)^2} \sum_{k=0}^\infty (-1)^k \gamma^k g_2((\Phi_1\Phi_2)^k \Phi_1 (x))\\
&\quad+ \frac{\beta}{k_2}  \log|x+(a,0)| \int_{\cB_2} \tilde f_3(y) \, dy\\
&\quad- \frac{2\alpha \beta}{k_2+1} \sum_{k=0}^\infty (-1)^k \gamma^k \log |(\Phi_1\Phi_2)^k \Phi_1(x) + (a,0)|\int_{\cB_2} \tilde f_3(y) \, dy\\
&\quad+ \frac{2}{k_2+1} \sum_{k=0}^\infty (-1)^k \gamma^k \Big( g_0((\Phi_2\Phi_1)^k (x)) - \alpha g_0((\Phi_1\Phi_2)^k \Phi_1 (x)) \Big),
\end{align*}
where we used $\log |x+(a,0)| = - \log |\Phi_2(x) + (a,0)|$ for $x \in \cB_2 \cap B_{1/2}$.
Note that
\begin{align*}
|x+(a,0)| \ge 1/2, \quad |(\Phi_1\Phi_2)^k \Phi_1(x) + (a,0)| \ge 1, \quad |(\Phi_2\Phi_1)^k (x)-(a,0)| \ge 1.
\end{align*}
We can estimate $Dw_i,i=0,1,2,3,$ as above to obtain
\begin{align*}
&|Dw_1(x)|\le \frac{Ck_1}{(k_1+1)(k_2+1)},\quad Dw_2(x)| \le C,\\
&|Dw_0(x)| \le \frac{C}{k_2+1},\quad |Dw_3(x)| \le C.
\end{align*}
Therefore,
$$
|Du(x)| \le C \quad \mbox{for}~~x \in \cB_2 \cap B_{1/2},
$$
where $C$ is a positive constant independent of $\varepsilon$, $k_1$ and $k_2$. Hence \eqref{main_estimate} is proved for $m=1$.

\section{Higher order derivatives}

In this section, we prove Theorem \ref{main_thm} for $m \ge 2$. The idea is essentially the same as the case when $m = 1$. Without loss of generality, we assume that $C_m = 1$.

For $x \in \cB_0 \cap B_{1/2}$, by \eqref{green_function} and \eqref{tilde_u_def}, we have
\begin{align*}
D^m w_1(x) = \frac{2}{k_1+1} \sum_{k=0}^\infty (-1)^k \gamma^k \Big( D^m[h_1(\psi^k (x))] - \beta D^m[h_1(\psi^k \Phi_2 (x))] \Big).
\end{align*}
We denote
$$\Theta_{m,k} := \gamma^k D^m[h_1(\psi^k (x))].$$
Since $\sum_{i=0}^m \binom{m}{i} = 2^m$, we have
\begin{align}\label{finite_difference}
\sum_{k=0}^{\infty} (-1)^k \Theta_{m,k} =& \frac{1}{2^m} \sum_{k=1}^\infty (-1)^k \left[ \sum_{i=0}^m (-1)^i \binom{m}{i} \Theta_{m, k+i} \right] \nonumber\\
&+ \Theta_{m,0} + \frac{1}{2^m} \sum_{i=1}^m \left[ (-1)^i \Theta_{m,i} \sum_{j = i}^m \binom{m}{j} \right].
\end{align}
The second line of \eqref{finite_difference} is clearly bounded independent of $\varepsilon$ and $\gamma$. The main goal is to estimate the first term on the right-hand side of \eqref{finite_difference}. By the mean value theorem,
$$
\sum_{i=0}^m (-1)^i \binom{m}{i} \Theta_{m, k+i} =  D_k^m \Theta_{m,k}\Big|_{k = \bar{k}}
$$
for some $\bar{k} \in (k, k+m)$. We claim that
\begin{equation}\label{goal}
\left| \sum_{k = 1}^\infty \sum_{i=0}^m (-1)^i \binom{m}{i} \Theta_{m, k+i} \right| \lesssim \sum_{k = 1}^\infty |D_k^m \Theta_{m,k}| \lesssim k_1.
\end{equation}
By the chain rule and the product rule,
$$
D^m[h_1(\psi^k (x))] = \sum_{n=1}^m \left( D^n h_1 (\psi^k(x)) \sum_{\sum_{j=1}^n \alpha_{n,j} = m} \prod_{j=1}^n D^{\alpha_{n,j}}(\psi^{k})(x) \right),
$$
where $\alpha_{n,j} \ge 1$.
For each $n$ and $\{\alpha_{n,j}\}_{j=1}^n$ satisfying $\sum_{j = 1}^n \alpha_{n,j} = m$,
$$
D^{\alpha_{n,j}}(\psi^{k})(x) = \sum_{\tilde\alpha_{n,j}=1}^{\alpha_{n,j}} \left( D^{\tilde\alpha_{n,j}}  \psi^{k-1}[\psi(x)] \sum_{\sum_{p=1}^{\tilde\alpha_{n,j}} \alpha_{n,j,p} = {\alpha_{n,j}}} \prod_{p=1}^{\tilde\alpha_{n,j}} D^{\alpha_{n,j,p}}\psi(x) \right),
$$
where $\alpha_{n,j,p} \ge 1$. Then by \eqref{D_psi},
\begin{align}\label{dkdx_expression}
& \Big| D_k^m \gamma^k D^m[h_1(\psi^k (x))] \Big|\nonumber\\
& \lesssim \sum_{\substack{\tau + \beta_0 = 0\\ \tau,\beta_0 \ge 0}}^m \sum_{n=1}^m \bigg( \gamma^k |\log \gamma|^\tau \big|D_k^{\beta_0} [D^n h_1 (\psi^k(x))]\big| \times \nonumber\\
&\quad \times \sum_{\sum_{j=1}^n \beta_j = m - \tau - \beta_0} \prod_{j=1}^n \sum_{\tilde\alpha_{n,j}=1}^{\alpha_{n,j}} \big|D_k^{\beta_j}D^{\tilde\alpha_{n,j}}(\psi^{k-1})[\psi(x)]\big| \bigg).
\end{align}
When $\beta_0 \ge 1$, by the chain rule and the product rule,
\begin{align*}
D_k^{\beta_0} [D^n h_1 (\psi^k(x))] = \sum_{s=1}^{\beta_0} \left( D^{n+s} h_1(\psi^k(x)) \sum_{\sum_{i=1}^s \gamma_{s,i} = \beta_0} \prod_{i=1}^s D_k^{\gamma_{s,i}} \psi^{k-1}[\psi(x)] \right),
\end{align*}
where $\gamma_{s,i} \ge 1$. For each $s$ and $\{\gamma_{s,i}\}_{i=1}^s$ satisfying $\gamma_{s,i} \ge 1$ and $\sum_{i = 1}^s \gamma_{s,i} = \beta_0$, by \eqref{dk}
\begin{align} \label{beta_0_positive}
\left|\prod_{i=1}^s D_k^{\gamma_{s,i}} \psi^{k-1}[\psi(x)] \right| &\lesssim \prod_{i=1}^s \left( \frac{\varepsilon^{(1 + \gamma_{s,i})/2}}{\lambda_2^{2k}} \sum_{j = 0}^{\gamma_{s,i}} \frac{r^{j+1}}{|I_{2k}|^{j+1}} \right)\notag\\
&\lesssim \frac{\varepsilon^{(s + \beta_0)/2}}{\lambda_2^{2k}} \sum_{j = 0}^{\beta_0} \frac{r^{j+s}}{|I_{2k}|^{j+s}}.
\end{align}
To prove the claim, we consider the following two cases:

\textbf{Case 1:} $\beta_0 \ge 1$. For each $n$, by \eqref{beta_0_positive}, \eqref{hg_derivatives}, and \eqref{dkdx}, the right-hand side of \eqref{dkdx_expression} can be estimated by (up to a positive constant independent of $k$)
\begin{align*}
&k_1\gamma^k |\log \gamma|^\tau \sum_{s=1}^{\beta_0} \left( \frac{\varepsilon^{(s + \beta_0)/2}}{\lambda_2^{2k}} \sum_{j = 0}^{\beta_0} \frac{r^{j+s}}{|I_{2k}|^{j+s}} \right) \frac{\varepsilon^{n + ( m - \tau -\beta_0)/2}}{\lambda_2^{2k}} \times\\
&\quad \times \sum_{\substack{0 \le \tilde\beta_2 + \tilde\beta_3 \le m - \tau - \beta_0\\
0 \le\tilde\beta_2 \le \tilde{m}-n}} |1- \lambda_2^{-2k}|^{\tilde{m}-n-\tilde{\beta}_2} \frac{r^{\tilde\beta_3}}{|I_{2k}|^{\tilde{m}+n+\tilde\beta_3}}\\
&\lesssim k_1\sum_{s=1}^{\beta_0} \gamma^k |\log \gamma|^\tau \frac{\varepsilon^{n + ( m - \tau + s)/2}}{\lambda_2^{2k}} \times\\
&\quad \times \sum_{\substack{ 0 \le\tilde\beta_2 + \tilde\beta_3 \le m - \tau - \beta_0\\
0 \le\tilde\beta_2 \le \tilde{m}-n\\
0 \le j \le \beta_0}} |1- \lambda_2^{-2k}|^{\tilde{m}-n-\tilde{\beta}_2} \frac{r^{\tilde\beta_3+j+s}}{|I_{2k}|^{\tilde{m}+n+\tilde\beta_3+j+s}},
\end{align*}
where $n \le \tilde{m} \le m$.
For every $s$, $\tau$, $j$, $\tilde{\beta}_2$, and $\tilde{\beta}_3$ that satisfy $1 \le s \le \beta_0$, $0 \le j \le \beta_0$, $0 \le\tilde\beta_2 \le \tilde{m}-n$, $\tau \ge 0$, and $0 \le\tilde\beta_2 + \tilde\beta_3 \le m - \tau - \beta_0$, we consider the following three cases:

When $k < \frac{1}{r+\sqrt{\varepsilon}}$, we have
$$
 |1- \lambda_2^{-2k}|\lesssim k \sqrt{\varepsilon} \quad \mbox{and} \quad |I_{2k}| \gtrsim \sqrt{\varepsilon}.
$$
Since
\begin{align*}
&\sum_{k < \frac{1}{r+\sqrt{\varepsilon}}} \left( \frac{\gamma}{\lambda_2^2} \right)^k k^{\tilde{m}-n- \tilde{\beta}_2} \\
&\lesssim \min \left\{ \frac{1}{(1 - \gamma/\lambda_2^2)^{\tilde{m}-n-\tilde{\beta}_2+1}}, \frac{1}{(r+ \sqrt{\varepsilon})^{\tilde{m}-n-\tilde{\beta}_2+1}} \right\}\\
&\lesssim \frac{1}{(1 - \gamma + r + \sqrt{\varepsilon})^{\tilde{m}-n- \tilde{\beta}_2+1}},
\end{align*}
we have
\begin{align*}
&\sum_{k < \frac{1}{r+\sqrt{\varepsilon}}} \gamma^k |\log \gamma|^\tau \frac{\varepsilon^{n + ( m - \tau + s)/2}}{\lambda_2^{2k}}  |1- \lambda_2^{-2k}|^{\tilde{m}-n-\tilde{\beta}_2} \frac{r^{\tilde\beta_3+j+s}}{|I_{2k}|^{\tilde{m}+n+\tilde\beta_3+j+s}}\\
&\lesssim \sum_{k < \frac{1}{r+\sqrt{\varepsilon}}}  \left( \frac{\gamma}{\lambda_2^2} \right)^k k^{\tilde{m}-n- \tilde{\beta}_2} |\log \gamma|^\tau \varepsilon^{(m + \tilde{m} + n - \tau - \tilde{\beta}_2 + s)/2} \frac{r^{\tilde\beta_3+j+s}}{\varepsilon^{(\tilde{m}+n+\tilde\beta_3+j+s)/2}}\\
&\lesssim \frac{1}{(1 - \gamma + r + \sqrt{\varepsilon})^{\tilde{m}-n- \tilde{\beta}_2+1}} (1 - \gamma)^\tau \varepsilon^{(m-\tau -\tilde{\beta}_2 - \tilde{\beta}_3-j)/2} r^{\tilde\beta_3+j+s}\\
&\lesssim \frac{1}{(1 - \gamma + r + \sqrt{\varepsilon})^{\tilde{m}-n- \tilde{\beta}_2+1}} (1 - \gamma + r + \sqrt{\varepsilon})^{m - \tilde{\beta}_2+s} \lesssim 1,
\end{align*}
as $\tilde{m} \le m$ and $n-1+s \ge 1 > 0$.
When $\frac{1}{r+ \sqrt{\varepsilon}} \le k < \frac{1}{\sqrt{\varepsilon}}$, we have
$$
 |1- \lambda_2^{-2k}| \lesssim k \sqrt{\varepsilon} \quad \mbox{and} \quad |I_{2k}| \gtrsim (r+\sqrt{\varepsilon})k\sqrt{\varepsilon}.
$$
Therefore,
\begin{align*}
&\sum_{\frac{1}{r+ \sqrt{\varepsilon}} \le k < \frac{1}{\sqrt{\varepsilon}}} \gamma^k |\log \gamma|^\tau \frac{\varepsilon^{n + ( m - \tau + s)/2}}{\lambda_2^{2k}}  |1- \lambda_2^{-2k}|^{\tilde{m}-n-\tilde{\beta}_2} \frac{r^{\tilde\beta_3+j+s}}{|I_{2k}|^{\tilde{m}+n+\tilde\beta_3+j+s}}\\
&\lesssim  \sum_{\frac{1}{r+ \sqrt{\varepsilon}} \le k < \frac{1}{\sqrt{\varepsilon}}} \gamma^{\frac{1}{r+ \sqrt{\varepsilon}}} k^{\tilde{m}-n- \tilde{\beta}_2} |\log \gamma|^\tau \varepsilon^{(m + \tilde{m} + n - \tau - \tilde{\beta}_2 + s)/2} \frac{r^{\tilde\beta_3+j+s}}{[(r+\sqrt{\varepsilon})k\sqrt{\varepsilon}]^{\tilde{m}+n+\tilde\beta_3+j+s}}\\
&\lesssim   \gamma^{\frac{1}{r+ \sqrt{\varepsilon}}} |\log \gamma|^\tau \varepsilon^{(m - \tau - \tilde{\beta}_2 - \tilde{\beta}_3 -j)/2} \frac{1}{(r+\sqrt{\varepsilon})^{\tilde{m}+n}} \sum_{\frac{1}{r+ \sqrt{\varepsilon}} \le k < \frac{1}{\sqrt{\varepsilon}}} \frac{1}{k^{2n+\tilde{\beta}_2 + \tilde{\beta}_3+s+j}}\\
&\lesssim  e^{\frac{\log \gamma}{r+\sqrt{\varepsilon}}} |\log \gamma|^\tau \varepsilon^{(m - \tau - \tilde{\beta}_2 - \tilde{\beta}_3 -j)/2} \frac{1}{(r+\sqrt{\varepsilon})^{\tilde{m}+n}} (r+\sqrt{\varepsilon})^{2n+\tilde{\beta}_2 + \tilde{\beta}_3+s+j-1} \\
&\lesssim  e^{\frac{\log \gamma}{r+\sqrt{\varepsilon}}} \left| \frac{\log \gamma}{r+\sqrt{\varepsilon}} \right|^\tau (r + \sqrt{\varepsilon})^{m - \tilde{m}+n+s-1} \lesssim 1,
\end{align*}
where we used the fact that $|e^{-x}x^\tau| \le C_\tau$ for any $x > 0$, $m - \tau - \tilde{\beta}_2 - \tilde{\beta}_3 -j \ge 0$, $\tilde{m} \le m$, and $n + s - 1 \ge 1 > 0$.

When $k \ge \frac{1}{\sqrt{\varepsilon}}$, we have
$$
 |1- \lambda_2^{-2k}| \lesssim 1 \quad \mbox{and} \quad |I_{2k}| \gtrsim r+\sqrt{\varepsilon}.
$$
Therefore,
\begin{align*}
&\sum_{k \ge \frac{1}{\sqrt{\varepsilon}}} \gamma^k |\log \gamma|^\tau \frac{\varepsilon^{n + ( m - \tau + s)/2}}{\lambda_2^{2k}}  |1- \lambda_2^{-2k}|^{\tilde{m}-n-\tilde{\beta}_2} \frac{r^{\tilde\beta_3+j+s}}{|I_{2k}|^{\tilde{m}+n+\tilde\beta_3+j+s}}\\
&\lesssim  \sum_{k \ge \frac{1}{\sqrt{\varepsilon}}} \left( \frac{\gamma}{\lambda_2^2} \right)^k |\log \gamma|^\tau\varepsilon^{n + ( m - \tau + s)/2} \frac{1}{(r+\sqrt{\varepsilon})^{m+n}}\\
&\lesssim \frac{\gamma^{\frac{1}{\sqrt{\varepsilon}}}}{1 - \gamma/\lambda_2^2} \left| \frac{\log \gamma}{\sqrt{\varepsilon}} \right|^\tau \varepsilon^{(n+s)/2}\\
&\lesssim e^{\frac{\log \gamma}{\sqrt{\varepsilon}}} \left| \frac{\log \gamma}{\sqrt{\varepsilon}} \right|^\tau \varepsilon^{(n+s-1)/2} \lesssim 1,
\end{align*}
where we used the fact that $|e^{-x}x^\tau| \le C_\tau$ for any $x > 0$, $\tilde{m}-n - \tilde{\beta}_2 \ge 0$, $n+s-1 \ge 1 > 0$, and $1 - \gamma/\lambda_2^2 \gtrsim \sqrt{\varepsilon}$. Therefore, \eqref{goal} follows from the estimates above.

\textbf{Case 2:} $\beta_0 =0$. For each $n$, by \eqref{dkdx}, the right-hand side of \eqref{dkdx_expression} can be estimated by (up to a positive constant independent of $k$)
\begin{align*}
k_1\gamma^k |\log \gamma|^\tau \frac{\varepsilon^{n + ( m - \tau)/2}}{\lambda_2^{2k}} \sum_{\substack{ 0 \le\tilde\beta_2 + \tilde\beta_3 \le m - \tau\\
0 \le\tilde\beta_2 \le \tilde{m}-n}} |1- \lambda_2^{-2k}|^{\tilde{m}-n-\tilde{\beta}_2} \frac{r^{\tilde\beta_3}}{|I_{2k}|^{\tilde{m}+n+\tilde\beta_3}},
\end{align*}
where $n \le \tilde{m} \le m$.
Then \eqref{goal} follows by repeating the exact same argument as in Case 1 with $s = 0$ and $j=0$.

Using the same argument but replacing $\psi^{k-1}(\psi(x))$ with $\psi^k \Phi_2(x)$, we have
$$
\left| \sum_{k=0}^\infty (-1)^k \gamma^k D^m[h_1(\psi^k \Phi_2 (x))] \right| \le Ck_1,
$$
where $C$ is a positive constant independent of $\varepsilon$ and $\gamma$. Therefore,
$$
|D^m w_1(x)| \le \frac{Ck_1}{k_1+1}\quad \text{for}\, x \in \cB_0 \cap B_{1/2}.
$$
As in Section 3, we can estimate $|D^m w_1(x)|$, $|D^m w_2(x)|$, $|D^m w_3(x)|$, and $|D^m w_0(x)|$ in all three regions. Therefore, \eqref{main_estimate} is proved.

\section{Proof of Theorem \ref{global_thm}}
In this section, we prove Theorem \ref{global_thm}. The proof is similar to that of Theorem \ref{main_thm}. Without loss of generality, we assume $C_m = 1$. Take  a domain $\cD_2$ such that $\cD_1 \Subset \cD_2 \Subset \cD$, and take a cutoff function $\eta \in C_0^\infty(\cD_2)$ such that $\eta = 1$ on $\cD_1$. Then $v := u \eta$ satisfies
$$
D_i(a(x) D_i v) = D_i \tilde f_i + \tilde{f}_3 \quad \mbox{in}~~\bR^2,
$$
where
\begin{equation*}
\tilde f_i = f_i \eta + u D_i \eta, \quad \tilde{f}_3 = - f_i D_i \eta + D_i u D_i \eta.
\end{equation*}
For $i=0,1,2,3$, $j = 1,2,3$, we define $\tilde{u}$, $w_i$, $h_j$, and $g_j$ as in \eqref{tilde_u_def}, \eqref{h_def}, and \eqref{g_def}. Instead of \eqref{uf_derivatives} and \eqref{hg_derivatives}, we have
\begin{equation*}
\| \tilde{f}_i \|_{C^{2m-1,\mu}(\cB_j)}\le C \min\{1,k_j\} \quad \mbox{for}~~i=1,2,3,~~ j=0,1,2,
\end{equation*}
and
\begin{align}
\|h_j\|_{C^{2m,\mu}(B_3 \cap \cB_i)} + \|g_j\|_{C^{2m,\mu}(B_3 \cap \cB_i)}\le C \min\{1,k_j\} \quad \mbox{for}~~i,j = 0,1,2 .\label{hg_derivatives_2}
\end{align}
As in the proof of Theorem \ref{main_thm}, for $x \in \cB_0 \cap \cD_1$, we have, by \eqref{hg_derivatives_2},
\begin{align*}
|D^m w_1(x)| \le \frac{Ck_1}{k_1+1}, \quad |D^m w_2(x)| \le \frac{C}{k_2+1}, \quad \mbox{and}~~ |D^m w_0(x)| + |D^m w_3(x)| \le C.
\end{align*}
For $x \in \cB_1$, we have
\begin{align*}
|D^m w_1(x)| &\le C + \frac{Ck_1}{(k_1+1)^2}, \\
|D^m w_2(x)| &\le \frac{C}{(k_1+1)(k_2+1)}, \\
|D^m w_0(x)| &\le \frac{C}{k_1+1}.
\end{align*}
To estimate $|D^m w_3(x)|$, we note that
\begin{align*}
&\frac{\alpha}{k_1} g_1(\Phi_1(x)) + \frac{\alpha}{k_1}\log|x-(a,0)| \int_{\cB_1} \tilde f_3(y) \, dy\\
&= \frac{\alpha}{k_1} \int_{\cB_1} (\log|\Phi_1(x) - y| + \log|x-(a,0)|) \tilde f_3(y) \, dy,
\end{align*}
and for fixed $y \in \cB_1$, $x = (a,0)$ is a removable singular point for $\log|\Phi_1(x) - y| + \log|x-(a,0)|$. Therefore $\log|\Phi_1(x) - y| + \log|x-(a,0)|$ is harmonic in $\cB_1$, and
$$
\left| D^m \left( \frac{\alpha}{k_1} g_1(\Phi_1(x)) + \frac{\alpha}{k_1}\log|x-(a,0)| \int_{\cB_1} \tilde f_3(y) \, dy \right) \right| \le C.
$$
The rest of the terms can be estimated similarly as before, and hence we obtain
$$
|D^m w_3(x)| \le C,
$$
which yields
$$
|D\tilde{u}(x)| \le \frac{C}{k_1 + 1} \quad \mbox{for}~~x \in \cB_1.
$$
Finally, for $x \in \cB_2$, we have
\begin{align*}
|D^m w_1(x)| &\le \frac{Ck_1}{(k_1+1)(k_2+1)}, \\
|D^m w_2(x)| &\le \frac{C}{k_2} + \frac{C}{(k_2+1)^2}, \\
|D^m w_0(x)| &\le \frac{C}{k_2+1},\\
|D^m w_3(x)| &\le \frac{Ck_1}{(k_1+1)(k_2+1)} +\frac{Ck_1}{k_2+1} + \frac{C}{k_2}+ \frac{C}{(k_2+1)^2}+ \frac{C}{k_2 + 1},
\end{align*}
which yields
$$
|D\tilde{u}(x)| \le \frac{C}{k_2 + 1} \quad \mbox{for}~~x \in \cB_2.
$$
Theorem \ref{global_thm} is proved.

\section{Proof of Theorem \ref{general_thm}}
When $r_1 = r_2$, Theorem \ref{general_thm} follows from Theorem \ref{global_thm} after scaling. When $r_1 \neq r_2$, we will find a conformal map $T: \bC \to \bC$ that maps $\cB_1$ and $\cB_2$ to circles of the same radius. Without loss of generality, we may assume $r_2 > r_1$.

Let $T(z) = \frac{1}{z- z_0}$, where $z_0 \in \bC$ and $z_0 = z_1 + z_2 i$. It is well known that if $z_0 \not\in \overline{\cB_1 \cup  \cB_2}$, $T$ maps $\cB_1 \cup  \cB_2$ to two disks. After a direct computation, we know that $T$ maps $\cB_1$ to the disk of center $\frac{\bar{z}_0 - \varepsilon/2 - r_1}{r_1^2 - |z_0 - \varepsilon/2 - r_1|^2}$, radius $\frac{r_1}{|r_1^2 - |z_0 - \varepsilon/2 - r_1|^2|}$, and maps $\cB_2$ to the disk of center $\frac{\bar{z}_0 + \varepsilon/2 + r_2}{r_2^2 - |z_0 + \varepsilon/2 + r_2|^2}$, radius $\frac{r_2}{|r_2^2 - |z_0 + \varepsilon/2 + r_2|^2|}$. We only need to find $z_0 = z_1 + z_2 i$ such that
$$
r_1(r_2^2 - |z_0 + \varepsilon/2 + r_2|^2) = r_2(r_1^2 - |z_0 - \varepsilon/2 - r_1|^2).
$$
This is equivalent to
\begin{align*}
&\left( z_1 - \Big( \frac{\varepsilon}{2} \frac{r_1+r_2}{r_2 - r_1} + \frac{2r_1r_2}{r_2 - r_1} \Big) \right)^2 + z_2^2 \\
&= \Big(\frac{2r_1r_2}{r_2 - r_1} \Big)^2 + \frac{2\varepsilon r_1r_2(r_1+r_2)}{(r_2-r_1)^2} + \frac{\varepsilon^2}{4} \left[ \Big(\frac{r_1+r_2}{r_2-r_1}\Big)^2 -1\right].
\end{align*}
We take
\begin{align*}
z_1 &= \sqrt{\Big(\frac{2r_1r_2}{r_2 - r_1} \Big)^2 + \frac{2\varepsilon r_1r_2(r_1+r_2)}{(r_2-r_1)^2} + \frac{\varepsilon^2}{4} \left[ \Big(\frac{r_1+r_2}{r_2-r_1}\Big)^2 -1\right]} \\
&\quad + \frac{\varepsilon}{2}\frac{r_1+r_2}{r_2 - r_1} + \frac{2r_1r_2}{r_2 - r_1},\\z_2 &= 0.
\end{align*}
It is easy to see that $z_1 > \varepsilon/2 + 4r_1$. Therefore $z_0 \not\in \overline{\cB_1 \cup  \cB_2}$, and we can choose domains $\cD_1$ and $\cD$ such that $\cB_1 \cup \cB_2 \Subset \cD_1 \Subset \cD$, and $z_0 \not\in \overline{\cD}$. Hence $T$ is smooth in $\cD$ and is the desired conformal map.

\section{The extreme case}
In this section, we prove Theorem \ref{extreme_thm}. The proof essentially follows that of \cite{BLLY}*{Theorem 1.1}, with some modifications.

By considering $u-C$ instead of $u$, we may assume the constant $C = 0$ in \eqref{main_problem_narrow}. For any $0<t<s<1$, let $\eta \in C_c^\infty(\Omega_s)$ be a cutoff function such that $\eta = 1$ in $\Omega_t$ and $|D \eta| \le C (s-t)^{-1}$. Multiplying $u\eta^2$ on both sides of \eqref{main_problem_narrow} and integrating by parts, we have
$$
\int_{\Omega_1} a^{ij}D_j u D_i u \eta^2 + 2 a^{ij} D_j u D_i \eta u \eta = 0.
$$
By Young's inequality,
$$
\int_{\Omega_t} |D u|^2 \le \frac{C}{(s-t)^2} \int_{\Omega_s \setminus \Omega_t} u^2.
$$
Since $u = 0$ on $\Gamma_-$, by the Poincar\'e inequality in the $x_n$ direction, we have
$$
\int_{\Omega_s \setminus \Omega_t} u^2 \le C(\varepsilon+s^2)^2 \int_{\Omega_s \setminus \Omega_t} |D u|^2.
$$
Therefore,
\begin{equation} \label{iteration}
\int_{\Omega_t} |D u|^2 \le C_0 \left(\frac{\varepsilon+s^2}{s-t} \right)^2 \int_{\Omega_s \setminus \Omega_t} |D u|^2.
\end{equation}
Let $t_0 = r \in (\sqrt{\varepsilon},1/2)$ and $t_j = (1-jr)r$ for $j \in \bN$. Taking $s = t_j, t = t_{j+1}$ in \eqref{iteration}, we have
$$
\int_{\Omega_{t_{j+1}}} |D u|^2 \le 4C_0 \int_{\Omega_{t_{j}} \setminus \Omega_{t_{j+1}}} |D u|^2.
$$
Adding both sides by $4C_0\int_{\Omega_{t_{j+1}}} |D u|^2$ and dividing both sides by $1+ 4C_0$, we have
$$
\int_{\Omega_{t_{j+1}}} |D u|^2 \le \frac{4C_0}{1+4C_0} \int_{\Omega_{t_{j}}} |D u|^2.
$$
Let $k = \lfloor \frac{1}{2r} \rfloor$ and iterate the above inequality $k$ times. We have
\begin{equation} \label{grad_u_L2}
\int_{\Omega_{r/2}}|D u|^2 \le \left( \frac{4C_0}{1+4C_0} \right)^k \int_{\Omega_{r}}|D u|^2 \le C\mu^{\frac{1}{r}} \int_{\Omega_{1}}|u|^2,
\end{equation}
where $\mu \in (0,1)$ and $C$ are constants depending only on $n$, $\sigma$, $\|h_1\|_{C^2}$, and $\|h_2\|_{C^2}$. For any $\bar{x} = (\bar x', \bar x_n) \in \Omega_{1/2}$, let $R = \varepsilon + h_1(\bar{x}') - h_2(\bar{x}')$. We make the change of variables by setting
$$
\left\{
\begin{aligned}
&y' = x' - \bar{x}',\\
&y_n = 2R \left( \frac{x_n - h_2(x') + \varepsilon/2}{\varepsilon + h_1(x') - h_2(x')} - \frac{1}{2} \right)
\end{aligned}
\right.
$$
for $x = (x', x_n) \in \{ |x' - \bar x'| < R, -\varepsilon/2 + h_2(x') < x_n < \varepsilon/2 + h_1(x') \}$. This change of variables maps the domain above to $Q_R$, where
$$
Q_s := \{ (y', y_n) \in \bR^n ~\big|~ |y'|<s , |y_n|<s\}
$$
for $s > 0$. Let
$$
(b^{ij}(y)) = \frac{(\partial_x y) \big(a^{ij}(x(y))\big) (\partial_x y)^t}{\det (\partial_x y)},
$$
$\tilde{b}^{ij}(y) = {b}^{ij}(Ry)$, and $\tilde{u}(y) = u(Rx)$. Then $\tilde{u}$ satisfies
$$
\left\{
\begin{aligned}
-\partial_i(\tilde b^{ij}(y) \partial_j \tilde u(y)) &=0 \quad \mbox{in } Q_{1},\\
\tilde b^{nj}(y) \partial_j \tilde u(y) &= 0 \quad \mbox{on } \{y_n = 1\},\\
\tilde u(y) &= 0 \quad \mbox{on } \{y_n = -1\}.
\end{aligned}
\right.
$$
It is straightforward to verify that
$$
\frac{I}{C} \le \tilde{b} \le CI \quad \mbox{and} \quad \|\tilde{b}\|_{C^{m-1,\alpha}(Q_1)} \le C,
$$
where $C$ is a positive constant depending only on $n$, $\sigma$, $m$, $\alpha$, and $C_{m,\alpha}$. Then by the Schauder estimate, we have
$$
\max_{-1 \le y_n \le 1} |D^m \tilde{u}(0',y_n)| \le C \|D \tilde{u}\|_{L^2(Q_1)},
$$
which gives, after reversing the change of variables,
\begin{equation}\label{grad_m_u}
|D^m u(\bar x)| \le C R^{1-m-n/2}\|D u\|_{L^2(\Omega_{|\bar x'| + R})}\quad \forall \bar{x} \in \Omega_{1/2},
\end{equation}
where $C$ is a positive constant depending only on $n$, $\sigma$, $m$, $\alpha$, and $C_{m,\alpha}$. Then by using \eqref{grad_u_L2} and \eqref{grad_m_u} with $r = 2 \max( \sqrt{\varepsilon}, |x'| + R)$, we conclude the proof.

%\section*{Declarations}
%\subsection*{Data availability}
%Data sharing not applicable to this article as no datasets were generated or analysed during the current study.
%\subsection*{Conflicts of interests}
%The authors have no relevant financial or non-financial interests to disclose.

\bibliographystyle{amsplain}

\end{document}